\documentclass[twocolumn]{article}
\usepackage[utf8]{inputenc}
\usepackage[left=2cm, right=2cm, top=2cm]{geometry} 
\usepackage[small]{titlesec}
\usepackage{etoolbox}

\makeatletter
\patchcmd{\ttlh@hang}{\parindent\z@}{\parindent\z@\leavevmode}{}{}
\patchcmd{\ttlh@hang}{\noindent}{}{}{}
\makeatother

\newcommand{\minimize}[3]{\begin{array}{rl}
	{\underset{#2}{\textrm{minimize}}} & \begin{aligned}[t] 
		#1
	\end{aligned}  \\[15pt]
	\textrm{subject to} &
	\begin{aligned}[t] 
		#3
	\end{aligned}
\end{array}}

\newcommand{\tinitial}{t_0}
\newcommand{\tgoal}{t_g}

\newcommand{\Mod}[1]{\ \mathrm{mod}\ ( #1)}

\usepackage{mathptmx} 
\usepackage{times} 
\usepackage{amsmath} 
\usepackage{amssymb}  
\usepackage{algorithm}
\usepackage{algpseudocode}
\usepackage{color}
\usepackage{makecell}
\usepackage{subfig, float}
\usepackage{tikz,pgfplots}
\usepackage{epstopdf}
\usepackage{diagbox}
\usepackage{hyperref}
\usepgfplotslibrary{external} 

\makeatletter
\let\NAT@parse\undefined
\makeatother
\usepackage[sort&compress,numbers,square]{natbib}

\makeatletter
\def\maketag@@@#1{\hbox{\m@th\normalfont\normalsize#1}}
\makeatother

\newlength\figureheight 
\newlength\figurewidth  

\newtheorem{myex}{Example}

\DeclareMathAlphabet{\mathcalOld}{OMS}{cmsy}{m}{n}
\DeclareMathOperator*{\minimizer}{minimize}
\DeclareMathOperator*{\subjecttoo}{subject\:to}

\title{\LARGE \bf
	Improved Optimization of Motion Primitives for Motion Planning in State Lattices
}

\date{}
\author{Kristoffer Bergman, Oskar Ljungqvist and Daniel Axehill}

\begin{document}
	
	\algblock{ParFor}{EndParFor}
	\algnewcommand\algorithmicparfor{\textbf{parfor}}
	\algnewcommand\algorithmicpardo{\textbf{do}}
	\algnewcommand\algorithmicendparfor{\textbf{end\ parfor}}
	\algrenewtext{ParFor}[1]{\algorithmicparfor\ #1\ \algorithmicpardo}
	\algrenewtext{EndParFor}{\algorithmicendparfor}

	\maketitle
	\thispagestyle{empty}
	\pagestyle{empty}

		\textbf{\textit{Abstract ---}}\textbf{In this paper, we propose a framework for generating motion primitives for lattice-based motion planners automatically. Given a family of systems, the user only needs to specify which principle types of motions, which are here denoted maneuvers, that are relevant for the considered system family. Based on the selected maneuver types and a selected system instance, the algorithm not only automatically optimizes the motions connecting pre-defined boundary conditions, but also simultaneously optimizes the end-point boundary conditions as well. This significantly reduces the time consuming part of manually specifying all boundary value problems that should be solved, and no exhaustive search to generate feasible motions is required. In addition to handling static a priori known system parameters, the framework also allows for fast automatic re-optimization of motion primitives if the system parameters change while the system is in use, e.g, if the load significantly changes or a trailer with a new geometry is picked up by an autonomous truck. We also show in several numerical examples that the framework can enhance the performance of the motion planner in terms of total cost for the produced solution.}

\section{Introduction} \label{sec:intro}
Motion planning for dynamical systems has become an important area of research not the least due to the recent high interest in autonomous systems such as self-driving cars, unmanned aerial vehicles and robotic manipulators. In general, finding an optimal motion in real-time for a system with nonlinear dynamics in an environment with obstacles is a challenging problem. Hence, most motion planning algorithms that are being developed to solve this problem aim at finding suboptimal solutions~\cite{lavalle2006planning}. Popular motion planning algorithms for systems subject to differential constraints include probabilistic~\cite{kuwata2008motion,karaman2013sampling} and deterministic methods~\cite{pivtoraiko2009differentially} with different types of guarantees on the quality of the produced solution. 

\begin{figure}
	\centering
	\setlength\figureheight{0.26\textwidth}
	\setlength\figurewidth{0.26\textwidth}
	\input{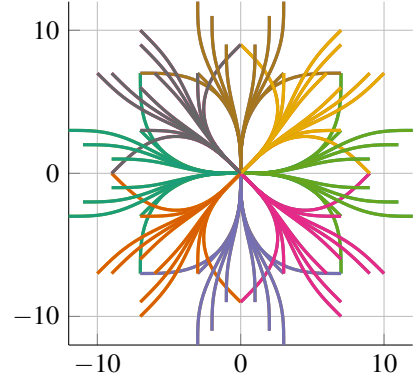}
	\caption{\label{fig:ex_prim} One example of a motion primitive set generated offline that can be used to solve motion planning problems online. }
\end{figure}

A popular deterministic method is the so called state lattice motion planning algorithm that was originally introduced in~\cite{pivtoraiko2009differentially} and has been used with great success on a variety of different vehicle platforms~\cite{ljungqvist2017lattice,BOSSDarpa,LjungqvistCDC2018}. The algorithm efficiently finds a solution to the motion planning problem by searching in a graph, where the vertices are discrete states and the edges are precomputed feasible motions that move the vehicle from one discrete state to another. These motions are called motion primitives and are generated offline by solving a set of Boundary Value Problems (BVPs). An example of a set of motion primitives is shown in Fig.~\ref{fig:ex_prim}. A solution to the motion planning problem is then given by an ordered sequence of motion primitives that connects the vehicle's current state to its desired goal state. The quality of the solution to the motion planning problem is thus determined by the discretization of the state space, the number of different motion primitives and the quality of each motion primitive~\cite{lindemann2005current}. 


In this paper, we consider parameterized families of nonlinear systems
\begin{equation}\label{eq:nonlinear_system}
\dot x(t) = f(x(t),u(t),p_f)
\end{equation}
where $x(t) \in \mathbb{R}^n$, $u(t) \in \mathbb{R}^m$, and $p_f \in \mathbb{R}^q$ denote the states, control signals, and the parameters that parameterizes the system family, respectively. Examples of parameters in common systems are platform dimensions or actuator limitations in cars, buses and trucks. The systems within the family are similar, but important parameters that affect optimal motions differs. If these variations are not handled in a structured way as proposed in this work, significant often manual work is necessary to develop state lattice motion planners for each instance of these system families since the motion primitives need to be recomputed using the correct values of the parameters. Using model-based numerical optimal control, this is fairly straight-forward given fixed boundary values. However, this work goes beyond that and introduces an efficient way to also optimize the end-point constraint itself. Intuitively, the optimal set of BVPs for a vehicle with a long wheel-base is not the same as for a vehicle with short wheel-base, since the former is far less agile compared to the latter. As will be shown in this work, for good performance it does not only suffice to optimize the behavior between the boundary conditions, but also the conditions themselves.

One approach to select the boundary conditions in the set of BVPs to be solved is by manual specification, which is typically done by an expert of the system~\cite{ljungqvist2017lattice, LjungqvistCDC2018}. The manual procedure can be very time-consuming and if some of the system's parameters are altered, the boundary conditions usually need to be carefully re-selected to not unnecessarily restrict the performance possible to obtain, and to maintain feasibility, during the following optimization of the motion primitives. Another widely used approach is to perform an exhaustive search to all neighboring states on the positional grid, either to the first feasible solution~\cite{pivtoraiko2009differentially} (shortest path set), or to all positions within a user-defined distance~\cite{pivtoraiko2005efficient,cirillo2014lattice}. These methods perform well for systems where analytical solutions exist, such as for differentially flat systems~\cite{sira2004differentially} where efficient exact steering methods exist~\cite{kelly2003reactive,murray1991steering}. However, for nonlinear systems that are not differentially flat, e.g., many truck and trailer systems~\cite{rouchon1993flatness}, the exhaustive search will lead to an intractable computational burden. 

To avoid time consuming manual tuning or exhaustive search-based generation of the motion primitive set for each system instance within the family, this work proposes a novel framework that is automatically optimizing both the boundary constraint selection in the BVPs and the corresponding motions to compute a motion primitive set for general nonlinear systems.  The inputs to the framework are an objective functional, the current system parameterization and a set of maneuvers that are specified by the user. It is shown in numerical examples that the framework significantly reduces the development time and enhances the generality of the motion primitive generation to allow for fast adaption to new system parameters. It is also illustrated that the optimized choice of connectivity allows the motion planner to generate motion plans with less total cost.

\section{Preliminaries}\label{sec:prel}
The optimal motion planning problem is to find a feasible trajectory that moves the vehicle from its initial state to a desired goal state, while not colliding with any obstacles and minimizes a specified performance measure. The problem can be formulated as a parametric continuous-time optimal control problem (OCP) according to 

\begin{equation}
\minimize{L(u(t), x(t), \tgoal, p_L)}{u(t), \;x(t), \;\tgoal}{&x(0) = x_{\mathrm{init}},
	\quad x(\tgoal) = x_{\mathrm{goal}}  \\ &\dot{x} (t) = f(x(t),u(t), p_f), \; \forall t \in [\tinitial,\tgoal] \\ &x(t) \in \mathcalOld{X}_{\mathrm{free}} \cap \mathcalOld{X}_{\mathrm{valid}}(p_x) , \; \forall t \in [\tinitial,\tgoal] \\
	& u(t) \in \mathcalOld{U}(p_u), \; \forall t \in [\tinitial,\tgoal]} \label{eq:cctoc}
\end{equation}
where $x_{\mathrm{init}}$ and $x_{\mathrm{goal}}$ are the initial and goal states, $\tgoal$ the time to reach the goal state, $\mathcalOld{X}_{\mathrm{free}}$ describes the obstacle free region, while $\mathcalOld{X}_{\mathrm{valid}}(p_x)$ and $\mathcalOld{U}(p_u)$ describe the physical constraints on the states and controls, and $L(u(t), x(t), \tgoal, p_L)$ forms the objective functional. The parameters $p = (p_L, p_x, p_u) \in \mathbb{R}^p$ in the problem represent parameters to the objective function, the state constraints and control input constraints, respectively,~\cite{buskens2000sqp}. In this paper, it is assumed that these parameters, together with the parameterization of the system $p_f$, remain constant during online planning. However, they are used in the problem formulation to specify the instance within a family of systems.

For fixed parameter values, the continuous-time OCP \eqref{eq:cctoc} is today intractable to solve online in real-time for generating long plans in environments with many obstacles. It is also necessary to be able to replan for example if the information about the surrounding environment is constantly updated~\cite{LjungqvistCDC2018}. Hence, approximate solutions in terms of discretized methods are commonly used~\cite{lavalle2006planning}. In this paper, a lattice-based approach is considered. 

\subsection{Lattice-based motion planning}
The idea with a lattice-based motion planner is to restrict the control to a discrete subset of the valid actions and hence transform the continuous-time motion planning problem to a discrete graph search problem. There are mainly two different approaches that are used to generate a lattice for motion planning~\cite{howard2008state}:
\begin{itemize}
	\item Control-sampling: The control space is sampled in a way such that the resulting sampling in the state space has desired properties in terms of discrepancy and dispersion. This typically lead to tree-shaped search spaces.   
	\item State lattice: First, a desired state space discretization $\mathcalOld{X}_d$ is selected. Then, a BVP solver is used to connect several states in the discretized state space. 
\end{itemize}
In this paper, the state lattice methodology will be used. The state lattice formulation has mainly been used for position invariant systems operating in unstructured environments. The benefits of using a state lattice are that it is possible to design the state space discretization depending on the application, and the complex relation between control and state dynamics is handled offline by the BVP solver. The latter also means reduced online computations since no forward simulations of the system are needed during online planning~\cite{pivtoraiko2011kinodynamic}. It is also much easier to create regular lattices, which will cover a larger volume of the search space in fewer samples~\cite{lavalle2006planning}. 
The use of a regular state lattice will lead to the possibility of using graph search methods for cyclic graphs, since many combinations of edges will arrive to the same state. For example, methods such as bidirectional search and exact pruning can be used~\cite{pivtoraiko2011kinodynamic}.

For a state lattice, the continuous-time OCP in \eqref{eq:cctoc} is transformed into a discrete graph search problem
\begin{subequations}
	\label{eq:planning}
	\begin{align}
	\minimizer_{\{m^k_p\}_0^{N-1},\hspace{0.5ex}N}\hspace{1ex}
	& \sum_{k=0}^{N-1} L_p(m_p^k)  \\
	\subjecttoo \hspace{1ex}
	& x_0 = x_{\mathrm{init}}, \quad x_N = x_{\mathrm{goal}}  \\ \label{eq:primStateTrans}
	& x_{k+1} = f_{m_p}(x_k,m_p^{u,k}), \; \forall k \in [0,N-1] \\ \label{eq:primInSet}
	& m_p^k  \in \mathcalOld{P}(x_k), \; \forall k \in [0,N-1] \\ \label{eq:trajInFreeSpace}
	& c(x_k,m_p^{x,k}) \in \mathcalOld{X}_{\mathrm{free}}, \; \forall k \in [0,N-1]
	\end{align} 
\end{subequations}   
where a motion primitive $m_p^k = (m_p^{x,k}, m_p^{u,k})$, $m_p^k \in \mathcalOld{P}$ is a dynamically feasible trajectory $\left( x^k(t), u^k(t)\right), t \in [ 0, T^k]$ that moves the vehicle from an initial state $x^k(0)\in \mathcalOld{X}_d$ to a final state $x^k(T^k)\in \mathcalOld{X}_d$, and the stage cost for each motion primitive is given by $L_p(m_p^k) = L(x^k(t), u^k(t), T^k, p_L)$. The equation in \eqref{eq:primStateTrans} defines the resulting state transition when $m_p^k$ is applied from $x_k$, and \eqref{eq:primInSet} ensures that $m_p^k$ is selected from the set of available motion primitives at $x_k$. Finally, the constraint in \eqref{eq:trajInFreeSpace} encodes that the resulting trajectory when $m_p^k$ is applied from $x_k$ has to remain in free-space.     

The motion primitive set $\mathcalOld{P}$ is generated offline. If the system is position invariant, it is only required to generate motion primitives from states with a position in the origin. These motion primitives can then be translated and reused at all positions in the grid~\cite{pivtoraiko2009differentially}. The construction of $\mathcalOld{P}$ is further described in the next section.

\section{State Lattice Construction}\label{sec:prim_gen}
\begin{figure}
	\centering
	\includegraphics[scale=0.6]{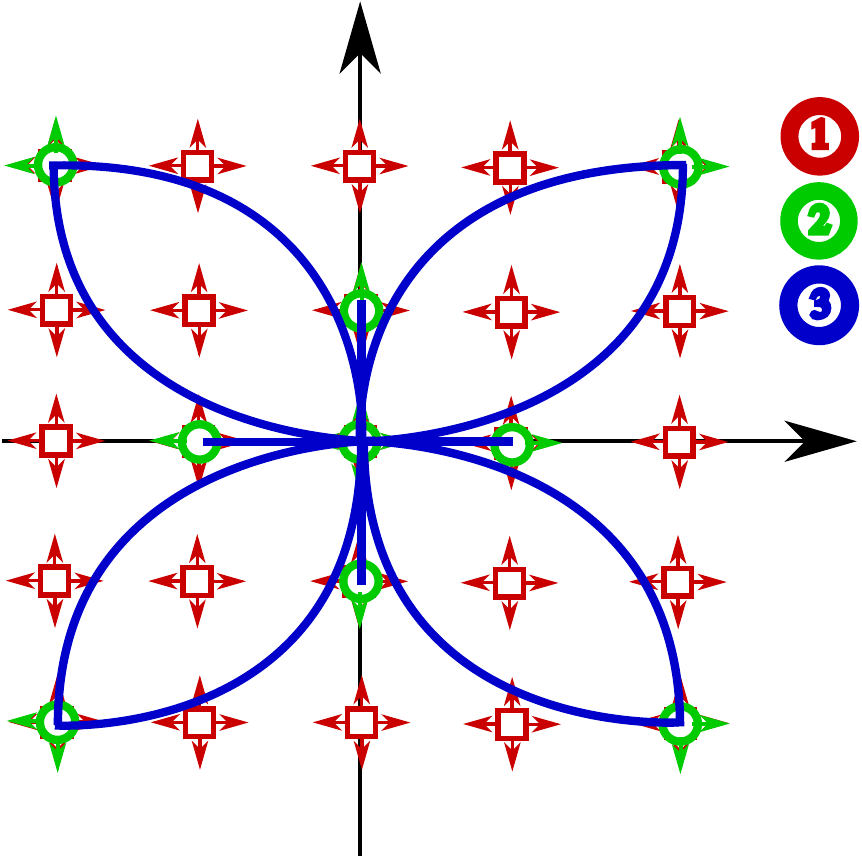}
	\caption{\label{fig:stepIll} An illustration of the three steps in the state lattice construction. (1) Discretize the state space, (2) select connectivity, (3) solve resulting BVPs.  }
\end{figure}

The aim of the state lattice construction is to compute a motion primitive set $ \mathcalOld{P}$ which contains feasible motions for the vehicle between different discretized vehicle states $x^d \in \mathcalOld{X}_d$ that can be used online to solve problems on the form \eqref{eq:planning}.
The construction procedure can be divided into three different steps~\cite{pivtoraiko2009differentially} (illustrated in Figure~\ref{fig:stepIll}):
\begin{enumerate}
	\item Discretize the state space to obtain a discrete search space,
	\item select which states to connect in the discretized representation,
	\item solve the BVPs defined in the previous step.
\end{enumerate}
The first step is to decide how the state space should be discretized, i.e., how to choose $\mathcalOld{X}^d$. This is done by selecting which fidelity of the state space that should be used in the search space for the planner. As an example for mobile robots, the first step could be to determine the fidelity of the position, orientation and steering angle of the vehicle~\cite{pivtoraiko2009differentially}. 

The second step is to select which states to connect in the discretized state space. Ideally, all possible initial states with a position in the origin to all final states in the discrete state space should be connected directly. However, this is intractable in practice since the required amount of memory usage, offline computation and online search time (with obstacles present) would be too high. Instead, a smaller number of connections are chosen, such that the resulting reachability graph $\mathcalOld{G}_r(x_{\mathrm{init}}, \mathcalOld{P})$ (which encodes the set of all possible trajectories from $x_{\mathrm{init}}$ given $\mathcalOld{P}$) sufficiently resembles the actual reachable set for the vehicle~\cite{lavalle2006planning}.

The procedure of how to select which states to connect is application dependent; systems that are more agile or systems working in environments with a lot of obstacles should have a dense connectivity selection. However,  it is also a trade-off between optimality and online computation, since the number of motion primitives defines the branching factor in the online graph search problem~\cite{pivtoraiko2009differentially}.

Finally, the BVPs defined by the second step needs to be solved to obtain the motion primitive set. In this paper, each motion primitive $m_p^i \in \mathcalOld{P}$ is defined as the solution to the following BVP

\begin{equation} 
\minimize{L(u^i(t), x^i(t), T^i, p_{L})}{u^i(t),x^i(t),T^i}{&x(0) = x^i_{\mathrm{init}}, \quad x^i(T^i) = x^i_{\mathrm{final}} \\ &\dot{x}^i (t) = f(x^i(t),u^i(t), p_f), \; \forall t \in [0,T^i] \\ &x^i(t) \in \mathcalOld{X}_\text{valid}(p_x) , \; \forall t \in [0,T^i] \\
	&u^i(t) \in \mathcalOld{U}(p_u), \; \forall t \in [0,T^i]\\
} \label{eq:primGen}
\end{equation}
where $x^i_{\mathrm{init}} \in \mathcalOld{X}^d$ and $x^i_{\mathrm{final}} \in \mathcalOld{X}^d$ are given by the connectivity selection in the previous step. Since these calculations are performed offline, the constraints representing the obstacles, $x(t) \in \mathcalOld{X}_{\mathrm{free}}$ in \eqref{eq:cctoc} are so far unknown and hence disregarded in \eqref{eq:primGen}. Instead, they are considered during online planning. The objective function can be chosen as any smooth function of the states, control signals and final time. For autonomous vehicles, this function is for example described as
\begin{equation} \label{eq:obj_ref}
L(u(t), x(t),T, \lambda, Q) = T + \lambda \int_{0}^{T} J\left(x(t),u(t), Q \right)  \mathrm{d}t
\end{equation}
where $J(x(t),u(t), Q) = ||\left[x(t), u(t) \right] ||_Q^2$ captures the smoothness behavior, and $ \lambda $ determines the trade-off between time duration and smoothness of a motion~\cite{ljungqvist2017lattice}.

One approach to solve the BVPs given in \eqref{eq:primGen} is to use numerical optimal control, which is a common method used for generating motion primitives~\cite{howard2007optimal,pivtoraiko2009differentially, ljungqvist2017lattice, LjungqvistCDC2018}. In this paper, we use state-of-the-art numerical optimal control software applicable for general nonlinear systems, such as ACADO~\cite{houska2011acado} or CasADi~\cite{andersson2018casadi}, to solve the optimal control problems. In these methods, the continuous-time optimal control problem in \eqref{eq:primGen} is reformulated into a nonlinear programming (NLP) problem using for example multiple shooting combined with numerical integration. A benefit of using numerical optimization is that state and control constraints can easily be incorporated. Also, it is straightforward to change the dynamical model of the system, the objective function, and to define and update problem parameters. Hence, it has the potential to be used as a backbone in a framework for generating motion primitives for a family of systems given the parameter that identifies the desired system instance.

\section{Maneuver-Based Motion Primitive Generation}\label{sec:autoGen}
In this section, the main contribution of this paper will be presented. The procedure of manually specifying the connectivity, i.e., step 2 in Section~\ref{sec:prim_gen}, is both time consuming and non-trivial since it is heavily application dependent. In this section, we propose a novel framework which simultaneously selects an optimized connectivity for a given system instance and computes optimal motions to generate the motion primitive set based on a number of user-defined maneuvers for which the state-connectivity is at last partially unspecified. 

\begin{figure}
	\centering
	\setlength\figureheight{0.3\textwidth}
	\setlength\figurewidth{0.3\textwidth}
	\input{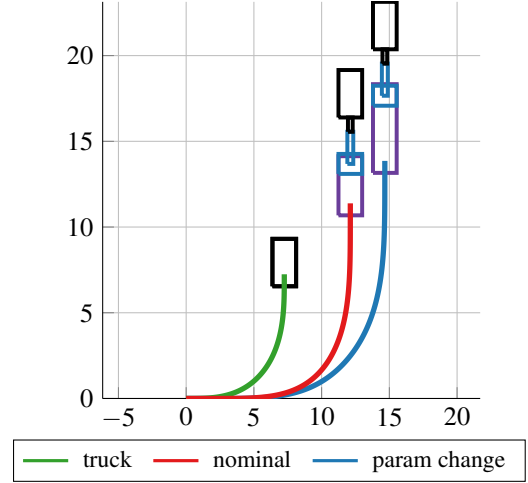}
	\definecolor{mycolor1}{rgb}{0.65098,0.80784,0.89020}%
\definecolor{mycolor2}{rgb}{0.12157,0.47059,0.70588}%
\definecolor{mycolor3}{rgb}{0.69804,0.87451,0.54118}%
\definecolor{mycolor4}{rgb}{0.20000,0.62745,0.17255}%
\definecolor{mycolor5}{rgb}{0.99216,0.74902,0.43529}%
\definecolor{mycolor6}{rgb}{0.892,0.102,0.1098}%
\definecolor{mycolor7}{rgb}{0.89020,0.10196,0.10980}%
\definecolor{mycolor8}{rgb}{0.41569,0.23922,0.60392}%

\begin{tikzpicture}

\begin{axis}[%
hide axis,
width=0,
height=0,
at={(0,0)},
scale only axis,
xmin=0,
xmax=1,
ymin=0,
ymax=1,
axis background/.style={fill=white},
xmajorgrids,
ymajorgrids,
legend style={legend cell align=left,column sep=3pt, align=center},
legend columns=3
]

\addplot [color=mycolor4, line width=1.2pt]
table[row sep=crcr]{%
	1	0\\
};
\addlegendentry{ {\small truck} }

\addplot [color=mycolor6, line width=1.2pt]
table[row sep=crcr]{%
	1	0\\
};
\addlegendentry{ {\small nominal} }
\addplot [color=mycolor2, line width=1.2pt]
table[row sep=crcr]{%
	1	0\\
};
\addlegendentry{ {\small param change} }

\end{axis}
\end{tikzpicture}%
	\caption{\label{fig:ex_comp} The optimal path, using the same objective function, for a $\pi/2$ heading change maneuver (without constraints on final position) changes depending on the system instance. From left to right: the optimal path for a truck, a truck with a dolly-steered trailer and finally the same truck and trailer system with a longer trailer. }
\end{figure}

Here, a \emph{maneuver} is defined as a principle type of desired motion for the vehicle system family that is parameterized. Formally, it is defined as a BVP similar to standard methods for computing motion primitives. The important difference is that in the BVP specifying the maneuver, freedom at the end-point has been introduced by removing one or several end-point constraints. This enables fast optimization to select an end-point constraint which is optimal for the selected system instance. The optimal choice of the end-point is application dependent; for example is the resulting optimal $\pi/2$ heading change maneuver for a truck is completely different from the solution for truck and trailer systems with different trailer lengths (illustrated in Fig.~\ref{fig:ex_comp}).  

From a maneuver specification and state space discretization, it is possible to easily compute a solution to any desired instance in the family. Depending on which type of maneuver, the final state constraint \mbox{$x(T)=x_{\mathrm{final}}$} in \eqref{eq:primGen} is replaced with a terminal manifold $ g(x(T)) = 0$ where $g:\mathbb{R}^n \rightarrow \mathbb{R}^l$ and $l<n$, where $n-l$ describes the degree of freedom for the terminal constraint, which is a part of the maneuver specification. This modification generalizes the BVP in \eqref{eq:primGen} for a traditional motion primitive to an OCP defining the more general maneuver
\begin{equation} 
\minimize{L(u^i(t), x^i(t), T^i, p_{L})}{u^i(t),x^i(t),T^i}{&x(0) = x^i_{\mathrm{init}}, \quad g^i(x^i(T^i)) = 0 \\ &\dot{x}^i (t) = f(x^i(t),u^i(t), p_f), \; \forall t \in [0,T^i] \\ &x^i(t) \in \mathcalOld{X}_\text{valid}(p_x) , \; \forall t \in [0,T^i] \\
	&u^i(t) \in \mathcalOld{U}(p_u), \; \forall t \in [0,T^i]\\
} \label{eq:primGenMan}
\end{equation}
To illustrate how a maneuver can be defined, assume that the vehicle under consideration can be described by the states $x = (x_1, y_1, \theta_1 )$ with $\dot{x} = \left[\cos(\theta _1), \sin (\theta _1), u_{\theta} \right]^T$. Assume that a discretized state space 

\begin{figure}
	\centering
	\setlength\figureheight{0.3\textwidth}
	\setlength\figurewidth{0.3\textwidth}
	\input{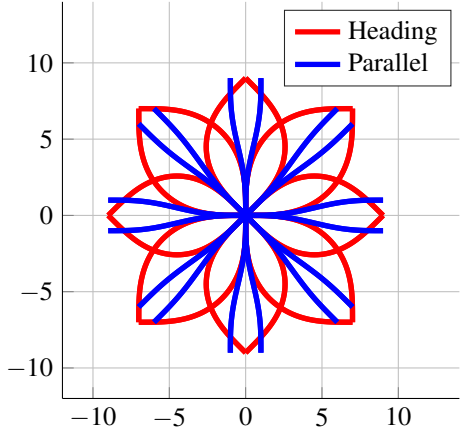}
	\caption{\label{fig:ex_man} An illustration of the resulting motion primitives from the maneuvers specified in Example 1-2 when $u_{\theta}$ is constrained, with $c_{\mathrm{lat}} = \pm 1$, $\Delta_{\theta} = 2$, $r=1$ and $\Theta = \{ \frac{k \pi}{4}\}_{k=0}^{7}$. }
\end{figure}

\begin{equation} \label{eq:ex_searchSpace}
x^d = \left(x_1^d, y_1^d, \theta_1^d \right)
\end{equation}
is given, where the position $(x_1^d, y_1^d)$ is discretized to a uniform grid with resolution $r$ and $\theta_1^d \in \Theta = \{ \theta_{1,k} \}_{k=0}^{N-1}$.     
\begin{myex} \label{ex:h_change}
	A \emph{Heading change maneuver} in \eqref{eq:ex_searchSpace} is given by a user-defined heading change to adjacent headings in $\Theta$. The boundary values are given by $x_{\mathrm{init}} = \left[0,0,\theta_{1,k} \right]^T$,  $\forall \; \theta_{1,k} \in \Theta$, to headings $\theta_{1,f} $, where $f=\Mod{k\pm\Delta_{\theta},N-1}$ and $\Delta_{\theta}$ is the user-defined parameter. The final state constraint then becomes
	\begin{equation}
	g(x(T)) = \theta_1(T) - \theta_{1,f} = 0,
	\end{equation} 
	while the values of $x_1(T)$ and $y_1(T)$ are left as free variables to be optimized by the OCP solver. The red paths in Fig.~\ref{fig:ex_man} show one example of motion primitives for a heading change maneuver.
\end{myex}

\begin{myex} \label{ex:p_man}
	A \emph{Parallel maneuver} in \eqref{eq:ex_searchSpace} is given by a user-defined lateral movement of $ c_{\mathrm{lat}}$ to the same heading as the initial one. The boundary values are thus given by $x_{\mathrm{init}} = \left[0,0,\theta_{1,k} \right]^T$,  $\forall \; \theta_{1,k} \in \Theta$, to  $\theta_{1,f} = \theta_{1,k}$, and with a final constraint that restrict  $(x_1(T), y_1(T) )$ to the line given by $\cos (\theta_{1,k})y_1(T) + \sin(\theta_{1,k})x_1(T) =  c_{\mathrm{lat}}$. Thus, the final state constraint for a parallel maneuver can be defined as
	\begin{equation}
	g(x(T)) = \begin{bmatrix}
	\theta_1(T) - \theta_{1,k} \\
	\cos (\theta_{1,k})y_1(T) + \sin(\theta_{1,k})x_1(T) - c_{\mathrm{lat}}
	\end{bmatrix} = 0.
	\end{equation}
	The blue paths in Fig.~\ref{fig:ex_man} show one example of motion primitives for a parallel maneuver.
\end{myex}

\begin{figure}[]
	\centering
	\includegraphics[scale=0.5]{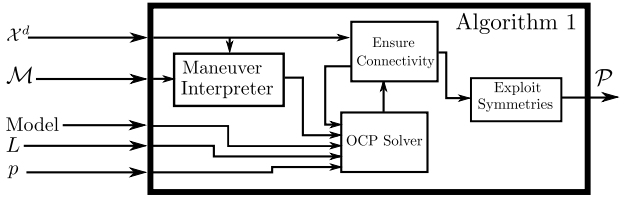}
	\caption{\label{fig:workflow} A flow chart of the proposed framework to generate motion primitives for a parameterized family of systems given a set of user-defined maneuvers.  }
\end{figure}

A flow chart of the proposed framework is shown in Fig.~\ref{fig:workflow} and described in Algorithm~\ref{alg:mbmpg}. The inputs are given by the state discretization, the system model and objective functional with parameters $p$ and a set of user-defined maneuvers $\mathcalOld{M}$. These inputs are used to setup the OCP solver as well as to define the initial and end-point constraints for each maneuver (line 2 and 3 in Algorithm~\ref{alg:mbmpg}), which results in a number of OCPs to be solved (one for each maneuver).   

One aspect to consider is that the resulting problem becomes a Mixed Integer Nonlinear Program (MINLP) due to the fact that formally the free final states are only allowed to take on the discrete values defined by the state discretization. To ensure that the motion primitives are feasible with the given discretization, a heuristic commonly used for solving MINLPs is used~\cite{bertsimas2005optimization}. A \emph{rounding heuristic} is applied from the obtained continuous OCP solution $x^i(t)$ where the closest end-point alternatives (in Euclidean sense) represented in $\mathcalOld{X}^d$ from $x^i(T^i)$ are evaluated and the result with lowest objective function value is chosen as resulting motion primitive $m_p^i$ (line 7 in Algorithm~\ref{alg:mbmpg}).

Finally, if the system is orientation invariant, rotational symmetries of the system can be exploited to reduce the number of OCPs that need to be solved, in the same way as presented in~\cite{pivtoraiko2005efficient} (here, line 8 in Algorithm~\ref{alg:mbmpg}). The symmetries of the system is something that should be accounted for already in the maneuver interpretation where the number of OCPs that need to be solved is decided.  

The benefits of using the proposed method outlined in Algorithm~\ref{alg:mbmpg} instead of completely manually specifying the connectivity are:
\begin{itemize}
	\item[+] The maneuver definitions can be \textbf{reused} to compute motion primitives for any desired instance of a parameterized family of systems, which reduces the manual workload of the operator since the re-connection of the state lattice will be handled automatically by the proposed framework.
	\item[+] Both the connectivity and the motions will be \textbf{optimized} for the specific system instance which gives a lower total cost function value in the online planning phase.
	\item[+] It is possible to automatically generate different levels of aggressiveness for the same maneuver automatically by simply changing $\lambda$ in \eqref{eq:obj_ref}, which can be useful to obtain smooth final solutions while still being able to express the reachability of the vehicle sufficiently well.
\end{itemize}

Note that the proposed framework goes beyond solving an OCP with new parameters, it also automatically optimizes the connectivity in the graph since some maneuver specifications, e.g. the end position, are left to be optimized by the framework. Industrially, motion planners can easily be prepared for new products with, e.g., another wheel-base. In that situation, the entire motion primitive set, including the state-connectivity, can be re-optimized only requiring a ``single-click'' from the developer.

\begin{algorithm}[t]
	\caption{Maneuver-Based Motion Primitive Generation}
	\label{alg:mbmpg}
	\begin{algorithmic}[1]
		\State \textbf{Input}: $\mathcalOld{X}^d$, system instance ($\mathcalOld{X}, \mathcalOld{U}, f$) and objective functional $L$ with parameters $p$ and a set of maneuvers $\mathcalOld{M}$ 
		\State OCP\_solver = setup\_problem($\mathcalOld{X}, \mathcalOld{U}, f, p, L$) 
		\State nOCP, $\mathbf{x}_{\mathrm{init}}, g(\mathbf{x}_{\mathrm{final}}) \leftarrow $ maneuver\_interpreter($\mathcalOld{M}$, $\mathcalOld{X}^d$  )
		\For{$i=1:\mathrm{nOCP}$}
		\State OCP\_solver.set\_boundary\_constraints$\left(x^i_{\mathrm{init}},g^i(x^i_{\mathrm{final}}) \right)$ 
		\State $x^i$ = OCP\_solver.solve()
		\State $m_p^i$ = ensure\_connectivity($x^i$, $\mathcalOld{X}^d$, OCP\_solver)
		\State $\mathcalOld{P}^i$ = exploit\_symmetries($m_p^i$)
		\State $\mathcalOld{P} = \mathcalOld{P} \cup \mathcalOld{P}^i$ 
		\EndFor	
	\end{algorithmic}
\end{algorithm}

Algorithm~\ref{alg:mbmpg} shares similarities with the motion primitive set generation frameworks described in \cite{pivtoraiko2009differentially, pivtoraiko2005efficient}. The method in~\cite{pivtoraiko2009differentially} generates a minimal set of motion primitives that connects all neighboring discrete states that do not represent the position of the vehicle by searching for the nearest feasible solution, and the method in~\cite{pivtoraiko2005efficient} searches for feasible end positions within a user-defined distance. Both are performed by an exhaustive search for feasible end positions by cycling through all possible candidates starting at the origin. However, the framework in this paper has some major advantages compared to the previously suggested methods:
\begin{itemize}
	\item[+] The computation time for the proposed method is orders of magnitude faster, and it scales better with state space discretization and problem dimension since the search for a feasible end-point in Algorithm~\ref{alg:mbmpg} starts at a dynamically feasible and optimal solution and not from the origin as in~\cite{pivtoraiko2009differentially, pivtoraiko2005efficient}. 
	\item[+]  Any user-defined (smooth) objective functional can be minimized. In \cite{pivtoraiko2009differentially}, the only objective function considered is minimum path length, since the search for a candidate end-point is terminated as soon as a feasible candidate is found. Therefore, the smoothness term in \eqref{eq:obj_ref} is always omitted in \cite{pivtoraiko2009differentially}.
\end{itemize}

\section{Numerical Examples} \label{sec:Res}
In this section, it is illustrated how the proposed framework can be applied to compute motion primitives for a diverse family of systems including cars and trucks with and without trailers. The implementation is done in Python, where CasADi combined with IPOPT~\cite{wachter2006implementation} is used as a backbone to solve the OCPs.
\subsection{Vehicle models}
The model for a car and a truck with no trailers is based on a kinematic car-like vehicle model~\cite{lavalle2006planning}:
\begin{equation} \label{eq:2dModel}
\small
\begin{aligned}
&\dot{x}_1 = v_1\cos \theta_1  \\
&\dot{y_1} = v_1\sin \theta_1   \\
&\dot{\theta_1}= v_1\tan \alpha  / L_1 \\
&\dot{\alpha} = \omega \\
&\dot{\omega} = u_{\alpha}
&\end{aligned}
\end{equation}
Here, $x = \left(x_1,y_1,\theta_1, \alpha, \omega \right)$ is the state vector which represents the position, heading, steering angle and steering angle rate, respectively. The control signals are $u = \left(v_1, u_{\alpha}\right)$, where $v_1$ is the longitudinal velocity of the rear axle of the vehicle and $u_\alpha$ the steering angle acceleration. Here, we consider path planning and thus constrain the longitudinal velocity as $v_1 = \{1, -1\}$ to generate both forward and backward motions. 

The truck and trailer system is a general 2-trailer with car-like truck that is illustrated in Fig.~\ref{fig:truck_model}. With the state vector $x_t = \left(x_3,y_3,\theta_3,\beta_3,\beta_2, \alpha, \omega \right)$, the kinematic model of this system is~\cite{altafini2002hybrid}:
\begin{figure}
	\centering
	\includegraphics[width=0.4\textwidth]{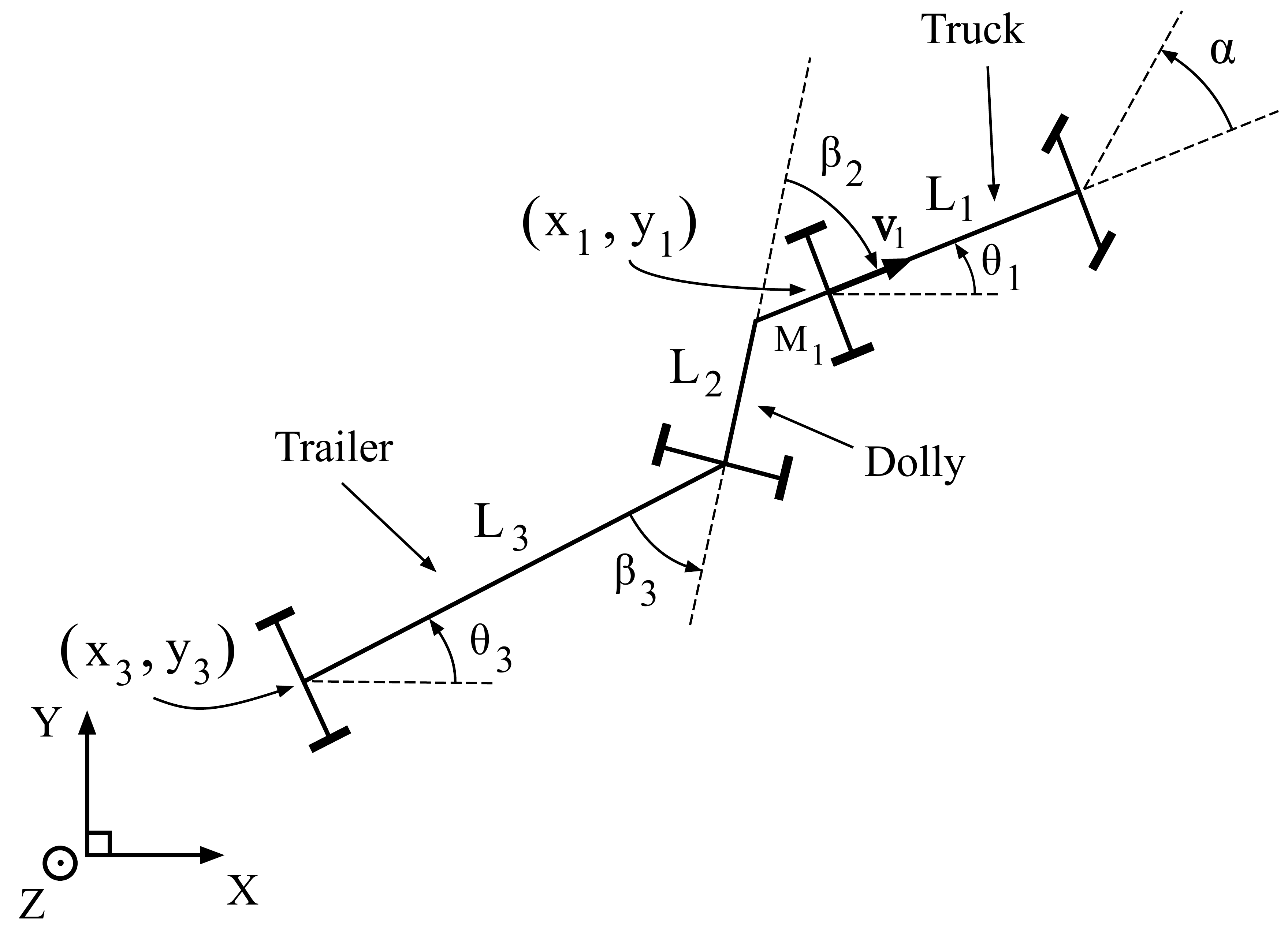}
	\caption{\label{fig:truck_model} A schematic illustration of the truck and trailer system that is used in the numerical examples. }
\end{figure}
\begin{equation} \label{eq:truckModel}
\small
\begin{aligned}
& \dot{x}_3 =  v_1 \cos \beta _3 \cos \beta _2 \left(1 + \frac{M1}{L1} \tan \beta _2  \tan \alpha \right) \cos \theta _3 \\ 
& \dot{y}_3 = v_1 \cos \beta _3 \cos \beta _2 \left(1 + \frac{M1}{L1} \tan \beta _2  \tan \alpha \right) \sin \theta _3 \\
& \dot{\theta}_3 = v_1 \frac{\sin \beta_3 \cos \beta _2}{L_3}\left(1 + \frac{M_1}{L_1} \tan \beta _2 \tan \alpha \right) \\
& \dot{\beta}_3 = v_1 \cos \beta_2 \Bigg(\frac{1}{L_2} \left( \tan \beta_2 - \frac{M_1}{L_1} \tan \alpha \right) - \\ 
& \quad \quad \quad  \frac{\sin \beta _3}{L_3} \left(1 + \frac{M_1}{L_1}\tan \beta_2 \tan \alpha \right) \Bigg) \\
& \dot{\beta}_2 = v_1 \left(\frac{\tan \alpha}{L_1} - \frac{\sin \beta_2}{L_2} + \frac{M_1}{L_1L_2} \cos \beta_2 \tan \alpha \right) \\ 
& \dot{\alpha} =  \omega \\
& \dot{\omega} = u_{\alpha}
\end{aligned}
\end{equation}
where $(x_3, y_3)$ and $\theta _3$ represent the center of the rear axle of the trailer and the heading of the trailer, respectively. Moreover, $\beta _3$ and $\beta_2$ denote the joint angles between the trailer and the truck. The control signals for the truck and trailer system coincide with the control signals for the car-like vehicle~\eqref{eq:2dModel}, \emph{i.e.}, $u = \left(v_1, u_{\alpha}\right)$. The system parameters for the general 2-trailer with car-like truck are the geometric lengths $L_3$, $L_2$, $L_1$ and $M_1$, which relate the distances and connections between the vehicle bodies. The truck and trailer system can be posed in \textit{circular equilibrium configurations}, where a constant steering angle $\alpha_e$ results in constant joint angles $\beta_{3,e}$ and $\beta_{2,e}$~\cite{altafini2002hybrid}. If such configurations are kept constant, the axles of the vehicle move along circles with radii determined by $\alpha_e$ and the geometric lengths $L_3$, $L_2$, $L_1$ and $M_1$. For more details regarding the model of the truck and trailer system, the reader is referred to~\cite{altafini2002hybrid}.     

The physical constraints on the steering angle, the steering angle rate and acceleration are the same for both systems, and are given by
\begin{equation*}
\mathcalOld{X}_{\mathrm{valid}}(p_x) = \left\{  x(t) :\begin{array}{lr}
-\alpha_{\mathrm{max}} \leq \alpha(t) \leq \alpha_{\mathrm{max}} \\
-\omega_{\mathrm{max}} \leq \omega(t) \leq \omega_{\mathrm{max}}
\end{array}\right\}
\end{equation*}
and $\mathcalOld{U}(p_u) = \{ u_{\alpha}(t) : -u_{\alpha,\mathrm{max}} \leq u_{\alpha}(t) \leq u_{\alpha,\mathrm{max}} \} $. In the numerical examples, the bounds are given by $ \alpha_{\mathrm{max}} = \pi/4$ $\mathrm{rad}$, $\omega_{\mathrm{max}} = 0.5$ $\mathrm{rad/s}$ and $u_{\alpha,\mathrm{max}} = 40$ $\mathrm{rad/s^2}$.
\subsection{State space discretizations}
In the state lattice formulation, the state space of the vehicles has to be discretized. The state space for the car-like vehicle is discretized as $x^d =  \left(x_1^d, y_1^d, \theta_1^d, \alpha^d, \omega^d \right)$, and for the truck and trailer system as $x^d_t =  \left(x_3^d, y_3^d, \theta_3^d, \beta_3^d, \beta_2^d, \alpha^d, \omega^d \right)$. The position for both vehicles $(x_i^d, y_i^d)$ are discretized to a uniform grid with resolution $r=1$. The headings are irregularly discretized $\theta_i^d \in \Theta = \{\theta_{k}\}_{k=0}^{N-1}$ to be able to generate short, straight motions from all headings~\cite{pivtoraiko2009differentially}. The steering angle is discretized according to $\alpha_e^d \in \Phi = \{\alpha_{e,k}\}_{k=0}^{M-1}$. For the truck and trailer case, in analogy with~\cite{ljungqvist2017lattice}, $\alpha_e^d$ is chosen to implicitly determine the discretization of the joint angles $\beta^d_{3,e}$ and $\beta^d_{2,e}$, where the vehicle is constrained to circular equilibrium configuration at each discrete state in the state lattice. This reduces the state dimension from seven to five states for the truck and trailer system, which imply that the state space discretization for both vehicles can be described by $x^d$. Finally, to ensure that $\alpha$ is a continuously differentiable function in time when multiple motions are combined during online planning, the steering angle rate $\omega^d$ is constrained to zero in the state space discretization.
\begin{figure}
	\centering
	\setlength\figureheight{0.25\textwidth}
	\setlength\figurewidth{0.25\textwidth}
	\input{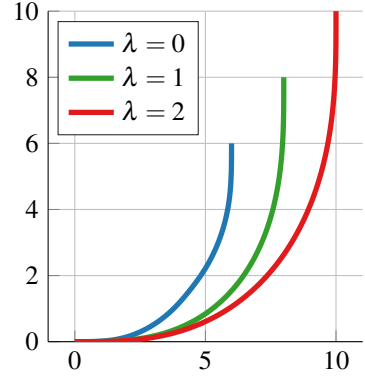}
	\caption{\label{fig:obj_comp} The resulting paths from an optimized heading change maneuver of $\pi/2$ for the truck using different values of $\lambda$ in the objective function \eqref{eq:objective}. Clearly, the choice of $\lambda$ has a significant impact on the resulting maneuver. }
\end{figure}

\subsection{Maneuver specifications}
The maneuvers that are typically used in a motion primitive set for the considered systems are 
straight, heading change, parallel or circular maneuvers. A straight maneuver is computed by connecting to the closest position on the grid that is aligned with the initial heading of the vehicle. The heading change and parallel maneuvers are defined in a similar way as in Example~\ref{ex:h_change}--\ref{ex:p_man}. The circular maneuvers are here used in combination with the heading change maneuver defined in Example~\ref{ex:h_change}, where each heading change maneuver is allowed to start from $\alpha^i_{e} \in \Phi$ and end in a neighboring $\alpha^f_{e} \in \Phi$. Since the state space discretization for the truck and trailer system is chosen such that $\alpha_e^d$ determines the values of $\beta_{3,e}^d$ and $\beta_{2,e}^d$, all maneuvers apply directly on both systems in \eqref{eq:2dModel} and \eqref{eq:truckModel}.  

The objective functional that is used in this section is given by
\begin{align} 
\small
L(x(t),u(t),T, \lambda) = T + \lambda \int_{0}^{T} \left( \alpha^2 + 10\omega^2 + u_{\alpha}^2 \right) \mathrm{d}t,
\label{eq:objective}
\end{align} 
When computing backward motions for the truck and trailer case, additional quadratic penalties for the joint angles $\beta_3$ and $\beta_2$ are added to promote solutions that have a small risk of leading to a jack-knife state during path execution, \emph{i.e}, 
\begin{equation} 
\small
\begin{aligned}
&L_b(x_t(t),u(t),T, \lambda) = \\
& T + \lambda \int_{0}^{T} \big( \beta_3^2 + \beta_2^2  + \alpha^2 + 10\omega^2 + u_{\alpha}^2 \big) \mathrm{d}t.
\label{eq:objectiveTruck}
\end{aligned}
\end{equation}
Further on in this section, $\lambda=1$ will be used unless stated otherwise. The impact of $\lambda$ on the optimal solution for the car-like vehicle can be seen in Fig.~\ref{fig:obj_comp} for a heading change maneuver of $\pi/2$. 
\subsection{Experimental results}
To illustrate the computational effectiveness of the proposed framework, four different sets of maneuvers are considered, which are defined and described in Table~\ref{tab:primSets}. The computation times for generating the complete set of motion primitives for several system instances are found in Table~\ref{tab:T_comp}. The system instances are a truck ($L_1=4.66$), a small car ($L_1=2.5$) and a truck and trailer system~\eqref{eq:truckModel} with two different trailer lengths ($L_3=8$ and $L_3=6$). For all settings using the proposed framework in this paper ($\mathcalOld{P}_1$-$\mathcalOld{P}_2$), the complete set of motion primitives are generated within minutes, which would have been significantly higher if the explicit state space connectivity for all maneuvers would have been chosen manually for each system instance. These results illustrate that the proposed framework could even be used online to enable automatic re-optimization of, e.g., an autonomous truck and trailer system if a new trailer with different length is to be connected to the truck. 


The proposed framework is also shown to outperform the suggested method in~\cite{pivtoraiko2009differentially} ($\mathcalOld{P}_3$-$\mathcalOld{P}_4$ in Table~\ref{tab:T_comp}), especially for the truck and trailer system. For this system instance, the algorithm in~\cite{pivtoraiko2009differentially} failed to compute the motion primitive set within ten hours. The main reason for the drastic change in computation times is that the method in~\cite{pivtoraiko2009differentially} spend most of the computations on trying to solve infeasible problems ($N_{\mathrm{inf}}$ in Table~\ref{tab:T_comp}) since it starts the search from the origin.  

\begin{figure}
	\centering
	\setlength\figureheight{0.25\textwidth}
	\setlength\figurewidth{0.35\textwidth}
%
%
\definecolor{mycolor1}{rgb}{0.41569,0.23922,0.60392}%
\definecolor{mycolor2}{rgb}{0.12157,0.47059,0.70588}%
\definecolor{mycolor3}{rgb}{0.89020,0.10196,0.10980}%
\definecolor{mycolor6}{rgb}{0.12157,0.47059,0.70588}%
\begin{tikzpicture}[arrow/.style={black,ultra thick,->,>=latex}]

\begin{axis}[%
width=\figurewidth,
height=\figureheight,
at={(0\figurewidth,0\figureheight)},
scale only axis,
xmajorgrids,
ymajorgrids,
xmin=-44.1653907825304,
xmax=30.7278907825304,
ymin=-35,
ymax=7.5,
axis background/.style={fill=white},
axis x line*=bottom,
axis y line*=left
]

\addplot [color=black, line width=1.5pt, forget plot]
  table[row sep=crcr]{%
-21.7484375	-23\\
-18.2515625	-23\\
-18.2515625	-17.46328125\\
-21.7484375	-17.46328125\\
-21.7484375	-23\\
};
\addplot [color=mycolor1, line width=1.5pt, forget plot]
  table[row sep=crcr]{%
-26.7484375	-35\\
-23.2515625	-35\\
-23.2515625	-26.95\\
-26.7484375	-26.95\\
-26.7484375	-35\\
};
\addplot [color=mycolor2, line width=1.5pt, forget plot]
  table[row sep=crcr]{%
-25.437109375	-28\\
-24.562890625	-28\\
-24.562890625	-24.25\\
-25.437109375	-24.25\\
-25.437109375	-28\\
};
\addplot [color=mycolor2, line width=1.5pt, forget plot]
  table[row sep=crcr]{%
-26.7484375	-29.165625\\
-23.2515625	-29.165625\\
-23.2515625	-26.834375\\
-26.7484375	-26.834375\\
-26.7484375	-29.165625\\
};
\addplot [color=black, line width=1.5pt, forget plot]
  table[row sep=crcr]{%
-26.7484375	-22.582\\
-23.2515625	-22.582\\
-23.2515625	-17.04528125\\
-26.7484375	-17.04528125\\
-26.7484375	-22.582\\
};
\addplot [color=black, line width=1.5pt, forget plot]
  table[row sep=crcr]{%
-25.3496875	-24.25\\
-24.6503125	-24.25\\
-24.6503125	-22.582\\
-25.3496875	-22.582\\
-25.3496875	-24.25\\
};

\addplot [color=black, line width=1.5pt, forget plot]
  table[row sep=crcr]{%
-16.125	-23\\
-12.875	-23\\
-12.875	-19.4375\\
-16.125	-19.4375\\
-16.125	-23\\
};

\addplot [color=black, line width=1.5pt, forget plot]
  table[row sep=crcr]{%
13	5.125\\
13	2.875\\
16.5625	2.875\\
16.5625	5.125\\
13	5.125\\
};
\addplot [color=black, fill=mycolor3,  forget plot]
  table[row sep=crcr]{%
-8	-1\\
0	-1\\
0	-5\\
-8	-5\\
-8	-1\\
};
\addplot [color=black, fill=mycolor3, forget plot]
  table[row sep=crcr]{%
-0.5	-13.5\\
4.5	-13.5\\
4.5	-23.5\\
-0.5	-23.5\\
-0.5	-13.5\\
};
\addplot [color=black, fill=mycolor3, forget plot]
  table[row sep=crcr]{%
10.5	-6.5\\
15.5	-6.5\\
15.5	-11.5\\
10.5	-11.5\\
10.5	-6.5\\
};
\addplot [color=black, fill=mycolor3, forget plot]
  table[row sep=crcr]{%
-13.5	7.5\\
-7.5	7.5\\
-7.5	4.5\\
-13.5	4.5\\
-13.5	7.5\\
};
\addplot [color=black, fill=mycolor3, forget plot]
  table[row sep=crcr]{%
6.5	4.5\\
9.5	4.5\\
9.5	1.5\\
12.5	1.5\\
12.5	-1.5\\
6.5	-1.5\\
6.5	4.5\\
};
\addplot [color=mycolor6,fill=mycolor6, fill opacity=0.1, dotted, line width=1.5pt, forget plot]
  table[row sep=crcr]{%
-30	-35.\\
-10	-35\\
-10	-16\\
-30	-16\\
-30	-35\\
};
\node[anchor=west] (source) at (axis cs:14,4){};
\node (destination) at (axis cs:23,4){};
\draw[->, ultra thick, >=latex](source)--(destination);
\node[] (source2) at (axis cs:-14.5,-20.5){};
\node (destination2) at (axis cs:-14.5,-14.5){};
\draw[->, ultra thick, >=latex](source2)--(destination2);
\node[] (source3) at (axis cs:-20,-18.5){};
\node (destination3) at (axis cs:-20,-12.5){};
\draw[->, ultra thick, >=latex] (source3) -- (destination3);
\node[] (source4) at (axis cs:-25,-18){};
\node (destination4) at (axis cs:-25,-12){};
\draw[->, ultra thick, >=latex] (source4) -- (destination4);
\end{axis}
\end{tikzpicture}%
	\caption{\label{fig:online_planning} Planning example with obstacles (red boxes) solved from several initial positions with $\theta_i=\pi/2$ (indicated by the area within the dotted lines) to a final state in $(x_i, y_i, \theta_i)= (13,4,0)$. The systems used in the example are (from left to right) a truck and trailer system, a truck ($L_1 = 4.66$) and a small car ($L_1=2.5$).   }
\end{figure}
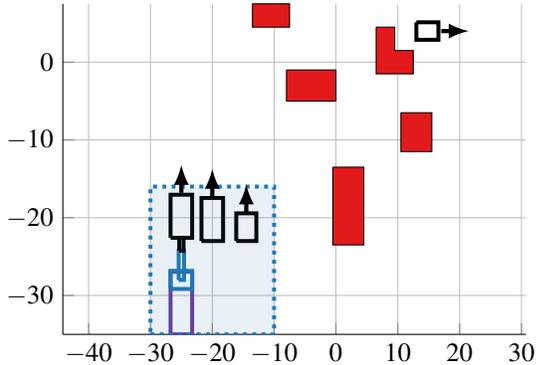

To demonstrate why it is important to choose the state space connectivity based on the platform to be controlled, a planning scenario from several initial states in an environment with obstacles is used (illustrated in in Fig.~\ref{fig:online_planning}) using a standard A$^*$ search to solve the motion planning problems. To ensure collision avoidance in the obtained solution, the obstacles are represented in the A$^*$ search using axis aligned bounding boxes, and the vehicles are represented by bounding circles~\cite{lavalle2006planning}. 

\begin{table}[t]
	\caption{A description of the different used motion primitive sets. $n_{\mathrm{\theta}}$ defines the number of $\theta^d$ in $\Theta$, $\Delta_{\theta}^{\mathrm{max}}$ defines which heading change maneuvers to generate (from 1 to $\Delta_{\theta}^{\mathrm{max}}$) and  $n_{\mathrm{par}}$ defines the number of parallel maneuvers from each initial heading. $\mathcalOld{P}_3$ and $\mathcalOld{P}_4$ contains no parallel maneuvers since they cannot be represented using the algorithm in~\cite{pivtoraiko2009differentially}.  Finally, $n_{\alpha_e}$ defines the number of $\alpha_e^d \in \Phi$ and $n_{\mathrm{prim}}$ the total number of motion primitives after exploiting system symmetries.} \label{tab:primSets}	
	\normalsize
	\centering
	\begin{tabular}{lllllllr}	
		\hline
		$\mathcalOld{P}$ & Alg. & Model &  $n_{\mathrm{\theta}}$ & $\Delta_{\theta}^{\mathrm{max}}$ & $n_{\mathrm{par}}$ & $n_{\alpha_e}$ & $n_{\mathrm{prim}}$ \\
		\hline
		$\mathcalOld{P}_1$ & Alg.~\ref{alg:mbmpg} & \eqref{eq:2dModel} & 16 & 4  & 4 & 3 & 1312   \\ 
		$\mathcalOld{P}_2$ & Alg.~\ref{alg:mbmpg} & \eqref{eq:truckModel} & 16 & 4 & 4  & 3 & 1312  \\ 
		$\mathcalOld{P}_3$ & \cite{pivtoraiko2009differentially} & \eqref{eq:2dModel} & 16 & 4 & 0 &  3 & 1056 \\ 
		$\mathcalOld{P}_4$   &  \cite{pivtoraiko2009differentially} & \eqref{eq:truckModel} & 16 & 4 & 0 & 3  & 1056 \\ 
		\hline
	\end{tabular}
\end{table}

In the comparison, the optimized connectivity for each one of the tested system instances are reused on each one of the remaining system instances. The results obtained from the planning scenario can be found in Table~\ref{tab:res}. When using a connectivity optimized for a system instance of smaller dimension, the motion primitive generation becomes infeasible, which is the case for the truck and trailer system and when using a connectivity optimized for a small car on a truck. In the opposite case, when using a connectivity optimized for a system instance of larger dimension, the motion primitive generation becomes feasible. However, in Table~\ref{tab:res} it can be seen that the average cost for the problems solved increases due to that the full potential of the maneuverability for the specific system instance is not fully utilized. Clearly, what is proposed in this work also has a value in practice.

\begin{table}[t]
	\caption{The total computation time ($T_{\mathrm{tot}}$) for generating and storing the complete set of motion primitives $\mathcalOld{P}_i$ defined in Table~\ref{tab:primSets}. The values of $L_1$, $L_2$, $L_3$ and $M_1$ represent the system parameter values used in \eqref{eq:2dModel} and \eqref{eq:truckModel}. $N_{\mathrm{inf}}$ represents the number of infeasible problems encountered during the motion primitive generation. } \label{tab:T_comp}	
	\normalsize
	\centering
	\begin{tabular}{llllllr}	
		\hline
		$\mathcalOld{P}$ & $L_1$ & $L_2$ & $L_3$ & $M_1$  & $T_{\mathrm{tot}}$ $\left[ s \right]$ & $N_{\mathrm{inf}}$ \\
		\hline
		$\mathcalOld{P}_1$ & 4.66 & N/A & N/A & N/A &  25.4 & 18  \\
		$\mathcalOld{P}_1$ & 2.5 & N/A & N/A & N/A &  28.6 & 36  \\	
		$\mathcalOld{P}_2$ & 4.66 & 3.75 & 8.0 & 1.67 & 314.2 & 208 \\
		$\mathcalOld{P}_2$ & 4.66 & 3.75 & 6.0 & 1.67 & 177.8 & 41 \\
		$\mathcalOld{P}_3$ & 4.66 & N/A & N/A & N/A  & 3947  & 21062 \\
		$\mathcalOld{P}_4$ & 4.66 & 3.75 & 6.0 & 1.67 & $>36000$ &   $> 46630$    \\  
		\hline
	\end{tabular}
\end{table}

\begin{table}
	\caption{Resulting relative cost for the planning scenario in Fig.~\ref{fig:online_planning} (comparable row-by-row). When a connectivity optimized for a system instance with larger dimensions is used, the vehicles are not able to utilize its maneuverability. In the opposite case, the motion primitive generation fails due to infeasibility. } \label{tab:res}
	\centering
	\normalsize
	\begin{tabular}{l*{3}{c}}\hline
		\backslashbox{Actual}{Opt. for}  
		&\makecell{$\mathcalOld{P}_2$ \\ $L_3=6$} & \makecell{$\mathcalOld{P}_1$ \\ $L_1=4.66$} & \makecell{$\mathcalOld{P}_1$ \\ $L_1=2.5$}   \\ \hline
		$\mathcalOld{P}_2$, $L_3=6$ & 1.0 & - & - \\
		$\mathcalOld{P}_1$, $L_1=4.66$ & 1.18  & 1.0  & - \\
		$\mathcalOld{P}_1$, $L_1=2.5$ & 1.31 & 1.14 & 1.0 \\ \hline
	\end{tabular}
\end{table}

\section{Conclusions and Future Work}
This work proposes a motion primitive generation framework for motion planning in state lattices. Based on user-defined principle motion types, here denoted maneuvers, the suggested framework automatically computes the motion primitive set for any user-selected system instance in a parameterized family of systems. This is performed by simultaneously optimizing the motion and the selection of the state-space connectivity. It is shown that this new framework enables the use of the same maneuver definitions for all instances in a parameterized, fairly diverse, family of systems, which enables fast automatized generation of the motion primitive set for any desired instance of the system. 

In numerical experiments the new method is shown to clearly outperform existing related methods, both in terms of performance as well as generality. The capabilities of the proposed framework both have applications at the level of industrial production of, e.g., several similar but different vehicles as well as online as a response to changes in the system to control. Furthermore, the proposed framework is shown in several numerical examples to increase the overall quality of the solutions generated in the online planning phase.

Future work includes to develop an extension to the proposed framework which automatically handles the discretization (completely or partially) of the state space, and to generate motion primitives that are also optimized with respect to expected potential obstacles in the environment. 

\section{Acknowledgments}
This work was partially supported by FFI/VINNOVA and the Wallenberg Artificial Intelligence, Autonomous Systems and Software Program (WASP) funded by Knut and Alice Wallenberg Foundation.
\bibliographystyle{abbrv}
\bibliography{myrefs.bib}		
\end{document}